\documentclass{article}

\usepackage[latin1]{inputenc}
\usepackage{amsmath, amsthm, amsfonts}
\usepackage[english]{babel}
\usepackage{graphics}
\usepackage{graphicx}
\usepackage{amsfonts}
\usepackage{amssymb}
\usepackage{subfigure}
\usepackage{ragged2e}

\theoremstyle{definition}

\theoremstyle{remark}

\newtheorem{defi}{Definition}[section]

\newtheorem*{Proof}{Proof}
\newtheorem*{teo2}{Keywords}

\newtheorem{lema}[defi]{Lemma}
\newtheorem{teo}[defi]{Theorem}

\theoremstyle{definition}

\theoremstyle{remark}
\newtheorem*{remark}{Remark}

\title{Spherical type surfaces via support function }
\author{Milton Javier Cardenas Mendez\and  Armardo Mauro Vasquez Corro}
\date{\today}

\begin{document}
\maketitle

\begin{abstract}
In this work we define the surfaces  spherical type via support function (in short, SS-surfaces).
We present a Weierstrass type representation for SS-surfaces with prescribed Gauss map which depends on two holomorphic functions. Also, we use this representation to classify the surfaces of rotation. Moreover, we show that every compact and connected SS-surface is the sphere and we give explicit examples of SS-surfaces.
\end{abstract}

\begin{teo2}
Generalized Weingarten surfaces,prescribed normal Gauss map, Weierstrass type representation.
\end{teo2}

\section{Introduction}
An oriented surface $\Sigma$ called a Weingarten surface if there is a differentiable relationship $W$ between the Gaussian curvature $K$ and the mean curvature H of S such that $W(H,K) = 0$.\\

The general classication of Weingarten surfaces is still an open question. When the functional W is linear,
namely, $ a + bH + cK = 0$ with constants $a, b, c \in \mathbb{R}$, the surfaces are called linear Weingarten surfaces.
Examples of linear Weingarten surfaces are the surfaces of constant Gaussian curvature ($b = 0$ and $c \neq 0$)
and the surfaces of constant mean curvature ($b  \neq 0$ and $c = 0$).\\

A surface $\Sigma$ is called the generalized Weingarten if it satisfies the following relation $A + BH + CK = 0$. Where $A$, $B$ and $C$ are functions that do not depend on the parameterization of $\Sigma$.\\

Let $\Sigma$ be a surface oriented by its normal Gauss map $N$. The functions $\psi,\Lambda :\Sigma\rightarrow \mathbb{R}^3$ given by
\begin{equation*}
\psi(p) = \langle p , N(p) \rangle \quad \text{and} \quad \Lambda(p) = \langle p , p \rangle \quad p \in  \Sigma,
\end{equation*}

where $\langle,\rangle$ denotes the Euclidean scalar product in $\mathbb{R}^3$, are called the support function and quadratic distance function, respectively.

Schief \cite{n} studied two classes of surfaces $ \Sigma \subseteq \mathbb{R}^3$ satisfying a Weingarten relation of the form $(\mu^2 ±\rho^2)K +
2\mu H+1 = 0$,where $\mu, \rho : \Sigma\rightarrow \mathbb{R}$  are harmonic functions in a certain sense. Bianchi \cite{g} classical surfaces, Bianchi
surfaces of positive curvature and the harmonic inverse mean curvature surfaces belong to these classes.\\

 In 1888 Appell \cite{f} studied a class of oriented surfaces in $\mathbb{R}^3$ associated with area preserving transformations in the sphere.Later, Ferreira and Roitman \cite{d} showed that these surfaces satisfy relation of Weingarten.
$ H +\psi K = 0$. Tzitzeica \cite{o} studied oriented hyperbolic surfaces such that there is a nonzero constant $c \in \mathbb{R}$ for which $K + c^2\psi^4 = 0$.

Martínez and Roitman in \cite{c}  show what seems to be the first example found for the second case of the problem proposed by Élie Cartan in their classic book a respect to external differential systems and their applications to Differential Geometry. Such examples are given by a class of Weingarten surfaces that satisfy the relationship $2\psi H+(1+\Lambda)K=0$

In \cite{a},the authors study a class of generalized special Weingarten surfaces, where coefficients are functions that depend on the support function and the distance function from a fixed point (in short EDSGW-surfaces), satisfy a relationship of the form $ A + BH + CK = 0 $, where $ A, B, C: S \rightarrow \mathbb{R} $ are differentiable functions.  He shows that, in the case where $ A, B, C $ are linear functions, a such a surface is invariant by inversions and expansions if, and only if, it satisfies a relation of the form $2\psi H +\Lambda K = 0$. In \cite{b}, the authors  present a Weierstrass type representation for EDSGW-surfaces with prescribed Gauss map which depends on two holomorphic functions. Also, we classify isothermic EDSGW-surfaces with respect to the third fundamental form parametrized by planar lines of curvature.

Let $ \Sigma $ be a hypersurface with $ \psi \neq 0 $, for each $ p \in \Sigma $, there is a sphere that touches $ \Sigma $ in $ p $ and passes through the origin.  This sphere has a radius given by $ R = - \frac{\Lambda} {2 \psi} $, which we will call the radius function.\\

Let $ \Sigma $ be a surface with normal application of Gauss $ N $, for each $ p \in \Sigma $ the center sphere $ p + \frac{H}{K}  N(p) $ and radius $ \frac{H}{K} $ is called \textbf {middle sphere}. A surface is called \textbf {spherical type} if the middle spheres tangent a fixed plane.\\

A surface is called the minimal Laguerre if $ \Delta_{III} (\frac{H}{K}) = 0 $. In 2009 Pottmann, Grohse and Mitra \cite{p}, proves that all spherical surfaces are minimal Laguerre.\\

Motivated by \cite{a},\cite{d},\cite{b}, we introduced a new class of surfaces. A $\Sigma $ surface is called a spherical type via support function (SS-surfaces) such that for each $ p \in \Sigma $ the center sphere $ p + (\frac{H(p)}{K (p )} + \frac{\psi(p)}{2}) N (p) $ and radius $ \frac{H(p)}{K(p)} + \frac{\psi(p)}{2 } $ go through the origin, in this case $ \Sigma $ satisfies the following generalized Weingarten relation
\begin{equation*}
2\psi H+(\Lambda+\psi^2)K=0
\end{equation*}
We characterize SS-surfaces in terms of a harmonic function and find a Weierstrass representation for these surfaces depending on two holomorphic functions and classify them in the case that they are rotational.

\section{Preliminaries}
In this section we fix the notation used in this work and present results of classical differential geometry of
surfaces. Throughout this paper $\Sigma$, $N$ and $U$ denote a surface in $\mathbb{R}^3$, its normal Gauss map and an open subset
of $\mathbb{R}^2$, respectively. The application $ X $ is called a parameterization of $ \Sigma $. If $ X_1, ..., X_{n + 1} $ indicate its component functions, then $ X $ is differential if and only if, $ X_i $ is differentiable for all $ i = 1, ..., n + 1 $ . In addition, the $ dX_q $ differential is injective if, and only if, the vectors
\begin{equation*}
X_{,i}(q) = \frac{\partial X} {\partial u_i}(q),1 \leq i \leq n
\end{equation*}
are linearly independent. \\
If $ X $ is a parameterization of $ \Sigma $, then $ T_p\Sigma $ matches the vector space generated by $ \{X_{,i}(u); 1 \leq i \leq n \}, \hspace {0.1cm} \text {where} \hspace{0.1cm} p = X(u) $. \\

Let $ X: U \subset \mathbb{R}^2 \rightarrow \Sigma,  $, a parameterization of a surface $ \Sigma $ and $ N: U \subset \mathbb{R}^2 \rightarrow \mathbb{R}^{3} $ a unitary vector field normal to $ X $, that is, $ N $ satisfies $ \langle N, X_{,i} \rangle (u) = 0 $, $ 1 \leq i \leq 2 $, for every $ u \in U $. Thus, the $ N $ field is normal to $ T_p\Sigma $, for each $ p = X (u) \in \Sigma $. We say that $ N $ is the normal Gaussian application of $ \Sigma $ and that such a field determines an orientation in $ \Sigma $.
Considering $ \{X_{,1}, X_{,2}, N \} $ as a base of $ \mathbb{R}^{3} $, we can write the
vector $ X_{,i j} $, $ 1 \leq i, j \leq 2 $, as
\begin{equation}
  X_{,ij}= \sum_ {k = 1}^ {2} \Gamma_ {ij}^k X_{,k} + b_{ij} N
\end{equation}

The $ \Gamma_{ij}^ k $ coefficients are called Christoffel symbols. If we take the parameterization $ X $ of $ \Sigma $ such that the metric $ L_ {ij} = \langle N_{,i}, N_{,j} \rangle $ is conform to Euclidean, that is, $ L_{ii } \neq 0 $ and $ L_{ij} = 0 $, for $ i \neq j $, The following lemma shows what Christoffel's symbols .

\begin{defi}
\label{def1}
Consider a $\Sigma$ surface of $ \mathbb{R}^{3} $ with map of Gauss $ N $. If $ X $ is a local parameterization of $ \Sigma $, the matrix $ W = (W_{ij}) $ such that
\begin{equation*}
N_{,i} = \sum_{j = 1}^{2} W_{ij} X_j, \hspace{0.3cm} 1 \leq i \leq 2
\end{equation*}
Is called the Weingarten matrix of $\Sigma$.
\end{defi}
Let $ Pr: \mathbb{S}^2 - \{(0, 0, -1)\}\rightarrow \Pi$ be the stereographic projection, starting from the south pole, of the sphere
$\mathbb{S}^2 = \{(x_1, x_2, x_3)\in \mathbb{R}^3 : x_1^2+ x_2^2+ x_3^2= 1\}$ in the plane $\Pi$. If $\Sigma$ has nonzero Gaussian curvature, then $X = Pr\circ N$ defines a local parametrization of plane $\Pi$. Thus, the normal Gauss map $N$ of $\Pi$ can be written as
\begin{equation*}
  N=\frac{( 2Y,1-\langle Y,Y\rangle)}{1+\langle Y,Y\rangle}
\end{equation*}
With $Y : U\rightarrow \Pi$  be an orthogonal local parametrization of the plane $\Pi$.
\begin{lema}
\label{le1}
Christoffel's symbols for the $ L_{ij} $ metric are given by
\begin{equation*}
\Gamma_{ij}^m=0 \hspace{0.3cm} \text{ for } i,j,k \hspace{0.3cm}  \text{distinct}
\end{equation*}
\begin{equation*}
\Gamma_{ij}^j=\frac{L_{jj,i}}{2L_{jj}} \hspace{0.3cm} \text{ for all } i,j
\end{equation*}
\begin{equation*}
\Gamma_{ii}^j=- \frac{L_{ii,j}}{2L_{jj}}=-\frac{L_{ii}}{L_{jj}}\Gamma_{ji}^i \hspace{0.3cm} \text{ for } i\neq j
\end{equation*}
\end{lema}

The next result obtained by Roitman and Ferreira  \cite{d}
\begin{teo}
\label{te1}
Let $ \Sigma \subset \mathbb{R}^{n + 1} $ be an orientable hypersurface and $ N: \Sigma \rightarrow \mathbb{S}^n $ normal Gauss application with non-zero Gauss-kroncker curvature in every point.  Let $ U \subset \Sigma $ be a neighborhood of $ p_0 $ such that $ N: U \rightarrow N(U) = V \subset \mathbb{S}^n $ invertible and $ h(q) = \langle q, N^{- 1}(q) \rangle, q \in V $ then
 \begin{equation}
\label{1}
X(q)=\nabla_Lh(q)+h(q)N(q)
\end{equation}
\end{teo}

\section{SS-surfaces}
Motivated by the works \cite{a},\cite{b},\cite{d}, we will start the study of Spherical type surfaces via support function and we call it SS-surfaces. In addition to presenting some examples, we provide a Weierstrass representation depending on two holomorphic functions for surfaces in this class and characterize the case where such surfaces are of rotation.

\begin{teo}
\label{te2}
Let $ h: U \subset \mathbb{R}^n \rightarrow \mathbb{R} $ be a differentiable function and $ N $ local orthogonal parameterization of $ \mathbb{S}^n $, then $ X $ given by (\ref {1}), defines a surface with normal application of Gauss N. Furthermore, the Weingarten matrix is given by
$ W = V^{- 1} $ where
\begin{equation}
\label{2}
V_{ij}=\frac{1}{L_{ij}}\left( h_{,ij}-\sum_{k}^{n}h_{,k}\Gamma_{ij}^{k}+hL_{ij}\delta_{ij}\right)
\end{equation}
And the fundamental forms $ I, II $ and $ III $ of $ X $, in local coordinates, are given by:
\begin{equation*}
 I=\langle X_{,i},X_{,j}\rangle=\sum_{k=1}^{n} V_{ik}V_{jk}L_{kk},\quad II=\langle X_{,i},N_{,j}\rangle=V_{ij}L_{jj},\quad III=\langle N_{,i},N_{,j}\rangle=L_{ij}\delta_{ij}
 \end{equation*}
\end{teo}

\begin{Proof}
Let $ X=\sum\limits_{j=1}^{n} \frac{h_{,j}}{L_{jj}}N_{,j}+hN $. Calculating the ith derivative and using  $\frac{L_{jj,i}}{L_{jj}}=2\Gamma_{ij}^j$, we obtain
\begin{equation*}
\begin{split}
  X_{,i}&=\left(\frac{h_{,ii}}{L_{ii}}-2\frac{h_{,i}}{L_{ii}}\Gamma_{ii}^i\right)N_{,i} + \frac{h_{,i}}{L_{ii}} \left(\sum_{k=1}^{n} \Gamma _{ii}^k N_{,k}-L_{ii}N \right)\\
 & +\sum_{\substack {j=1\\ j\neq i}}^{n}\left[\left(\frac{h_{,ji}}{L_{jj}}-2\frac{h_{,j}}{L_{jj}}\Gamma_{ij}^j\right)N_{,j} + \frac{h_{,j}}{L_{jj}} (\Gamma_{ji}^jN_{,j}+ \Gamma_{ji}^iN_{,i})\right]+h_{,i}N+hN_{,i}
\end{split}
\end{equation*}
As $\Gamma_{ij}^j=\Gamma_{ji}^j$ then
\begin{equation*}
\begin{split}
  X_{,i}&=\left(\frac{h_{,ii}}{L_{ii}}-2\frac{h_{,i}}{L_{ii}}\Gamma_{ii}^i\right)N_{,i} + \frac{h_{,i}}{L_{ii}} \sum_{k=1}^{n} \Gamma _{ii}^k N_{,k}- h_{,i}N  +\sum_{\substack {j=1\\ j\neq i}}^{n}\left(\frac{h_{,ji}}{L_{jj}}-\frac{h_{,j}}{L_{jj}}\Gamma_{ij}^j\right)N_{,j} -\sum_{\substack {j=1\\ j\neq i}}^{n}\frac{h_{,j}}{L_{jj}}\Gamma_{ij}^jN_{,j}\\
 & +\sum_{\substack {j=1\\ j\neq i}}^{n} \frac{h_{,j}}{L_{jj}} \Gamma_{ji}^jN_{,j}+\sum_{\substack {j=1\\ j\neq i}}^{n}\frac{h_{,j}}{L_{jj}} \Gamma_{ji}^iN_{,i}+h_{,i}N+hN_{,i}
\end{split}
\end{equation*}
rewriting
\begin{equation*}
X_{,i} =\left(\frac{h_{,ii}}{L_{ii}}-\frac{h_{,i}}{L_{ii}}\Gamma_{ii}^i+\sum_{\substack {j=1\\ j\neq i}}^{n}\frac{h_{,j}}{L_{jj}} \Gamma_{ji}^i +h \right)N_{,i} +\sum_{\substack {j=1\\ j\neq i}}^{n}\left(\frac{h_{,ji}}{L_{jj}}-\frac{h_{,j}}{L_{jj}}\Gamma_{ij}^j+\frac{h_{,i}}{L_{ii}}\Gamma_{ii}^j \right)N_{,j}
\end{equation*}
As $\Gamma_{ii}^j=-\frac{L_{ii}}{L_{jj}} \Gamma_{ji}^i$ for  $i\neq j$, we have to
\begin{equation*}
X_{,i}=\left(\frac{h_{,ii}}{L_{ii}}-\sum_{j=1}^{n}\frac{h_{,j}}{L_{ii}} \Gamma_{ii}^j +h \right)N_{,i}+\sum_{\substack {j=1\\ j\neq i}}^{n}\left(\frac{h_{,ij}}{L_{jj}}-\frac{h_{,j}}{L_{jj}}\Gamma_{ij}^j-\frac{h_{,i}}{L_{jj}}\Gamma_{ij}^i \right)N_{,j}
\end{equation*}
Considering the matrix $V=(V_{ij}),\hspace{0.3cm} 1\leq i,j\leq n$ given by
\begin{equation*}
V_{ij}=\frac{1}{L_{jj}}\left(h_{,ij}-\sum_{k=1}^{n}h_{,k}\Gamma_{ij}^k+hL_{ij}\delta_{ij}\right)
\end{equation*}
therefore
\begin{equation}
\label{14}
X_{,i}=\sum_{j=1}^{n}V_{ij}N_{,j}
\end{equation}
In search of the Weingarten matrix we have to $N_{,i}=\sum_{j=1}^{n}V^{-1}_{ij}X_{,j}$ by definition \ref{def1} we have to $W=V^{-1}$.
To obtain the coefficients of the fundamental forms, we use (\ref{14}), therefore
\begin{equation*}
I =\langle X_{,i},X_{,j}\rangle=\langle\sum_{k=1}^{n}V_{ik}N_{,k},\sum_{m=1}^{n}V_{jm}N_{,m}\rangle
=\sum_{k,m=1}^{n}V_{ik} V_{jm}\langle N_{,k}N_{,m}\rangle=\sum_{k=1}^{n}V_{ik} V_{jk}L_{kk}
\end{equation*}
\begin{equation*}
II =\langle X_{,i},N_{,j}\rangle=\langle\sum_{k=1}^{n}V_{ik}N_{,k},N_{,j}\rangle=\sum_{k=1}^{n}V_{ik}\langle N_{,k},N_{,j}\rangle=V_{ij}L_{jj}
\end{equation*}
\begin{equation*}
III =\langle N_{,i},N_{,j}\rangle=L_{ij}\delta_{ij}
\end{equation*}
\end{Proof}

\begin{defi}
A $\Sigma $ surface is called a spherical type via support function (SS-surfaces) such that for each $ p \in \Sigma $ the center sphere $ p + (\frac{H(p)}{K (p )} + \frac{\psi(p)}{2}) N (p) $ and radius $ \frac{H(p)}{K(p)} + \frac{\psi(p)}{2 } $ go through the origin, in this case $ \Sigma $ satisfies the following generalized Weingarten relation
\begin{equation*}
  2\psi H+(\Lambda+\psi^2)K=0
\end{equation*}
for all $p\in \Sigma$
\end{defi}
For SS-surfaces with Gaussian curvature $ K \neq 0 $ we will present a complete characterization through pairs of holomorphic functions. This representation will allow to classify all SS-surfaces of rotation. Before that, we will need the following lemma.

\begin{lema}
Consider holomorphic functions $g: \mathbb{C} \rightarrow \mathbb{C}_{\infty}$ and $f:\Sigma \rightarrow \mathbb{C}$, with $ g'\neq0 $, where $\Sigma$ is a Riemann surface. Taking the local parameters $ z = u_1 + iu_2 \in \Sigma $ e $ h = e^{\langle 1, f \rangle} $, the matrix
\begin{equation*}
V_{ij}=\frac{1}{L_{jj}}\left(h_{,ij}-\sum_{k=1}^{2}h_{,k}\Gamma_{ij}^k+hL_{ij}\delta_{ij}\right), \hspace{0.4cm} 1\leq i,j\leq 2
\end{equation*}
with
\begin{equation*}
L_{11}=L_{22}= \frac{4|g'|^2}{T^2},\quad L_{12}=L_{21}=0,\quad T=1+|g|^2
\end{equation*}
is such that
\begin{equation*}
\begin{split}
    V_{11} =\frac{T^2h}{4|g'|^2}[\langle 1,&f'\rangle^2-\langle 1,\xi \rangle] + h,\quad V_{12}=V_{21}=\frac{T^2h}{4|g'|^2}\langle i,\xi -\frac{f'^2}{2} \rangle  \\
     &   V_{22}=\frac{T^2h}{4|g'|^2}[\langle 1,if'\rangle^2+\langle 1,\xi \rangle] + h
\end{split}
\end{equation*}
Where $\xi =f' \left(\frac{g''}{g'}-\frac{2}{T}g'\overline{g}\right)-f''$, we have
\begin{equation}\label{3}
      trV=\frac{h|f'|^2T^2}{4|g'|^2} + 2h
\end{equation}
\begin{equation}\label{33}
  detV=\frac{T^4 h^2}{4|g'|^2}[\langle i,\frac{ f'^2}{2} \rangle^2-\langle i,\xi -\frac{f'^2}{2} \rangle^2+\langle 1,\xi \rangle\langle 1,f'^2-\xi \rangle+\frac{|f'|^2}{T^2}]+h^2
\end{equation}
\end{lema}

\begin{Proof}
With $h=e^{\langle 1,f\rangle}$, it follows
\begin{equation*}
\begin{split}
&h_{,1}=e^{\langle 1,f\rangle}\langle 1,f'\rangle, \quad h_{,11} = e^{\langle 1,f\rangle}(\langle 1,f'\rangle^2+\langle 1,f''\rangle),\quad h_{,2}=e^{\langle 1,f\rangle}\langle 1,if'\rangle,\\
     &h_{,22} = e^{\langle 1,f\rangle}(\langle 1,if'\rangle^2-\langle 1,f''\rangle),\qquad h_{,12} = e^{\langle 1,f\rangle}(\langle 1,if'\rangle\langle 1,f'\rangle+\langle 1,if''\rangle).
\end{split}
\end{equation*}
In order to find Christoffel's symbols for the metric, let's calculate the derivatives of the real functions  $T^2$ e $|g'|^2$. Since $T=1+|g|^2=1+\langle 1 ,g\rangle^2 +\langle i ,g\rangle^2$ we have
\begin{equation*}
\begin{split}
&T_{,1}=2\langle g,g'\rangle, \quad T_{,2}=2\langle g,ig'\rangle,\quad(T^2)_{,1}=4T\langle g,g'\rangle,\quad (T^2)_{,2}=4T\langle g,ig'\rangle.
\end{split}
\end{equation*}
On the other hand,
\begin{equation*}
|g'|_{,1} ^2 =(\langle 1,g'\rangle^2)_{,1}+ (\langle i,g'\rangle^2)_{,1} =2\langle 1,g'\rangle\langle 1,g''\rangle+2\langle i,g'\rangle\langle i,g''\rangle=2\langle g',g''\rangle
\end{equation*}
analogously   $|g'|_{,2} ^2 =2\langle g',ig''\rangle$.So, we have to
\begin{equation*}
   L_{11,1} =\left(\frac{4|g'|^2}{T^2}\right)_{,1} =\frac{8T\langle g',g''\rangle-16|g'|^2\langle g,g'\rangle}{T^3}
\end{equation*}
\begin{equation*}
   L_{22,2} =\left(\frac{4|g'|^2}{T^2}\right)_{,2} =\frac{8T\langle g',ig''\rangle-16|g'|^2\langle g,ig'\rangle}{T^3}
\end{equation*}
By the Lemma \ref{le1} we have  $\Gamma_{ij}^j=\frac{L_{jj,i}}{2L_{jj}}$  for all $i,j$
\begin{equation*}
   \Gamma_{11}^1 =\frac{L_{11,1}}{2L_{11}}=\frac{8T\langle g',g''\rangle-16|g'|^2\langle g,g'\rangle }{T^3} \frac{T^2}{8|g'|^2}=\frac{T\langle g',g''\rangle-2|g'|^2\langle g,g'\rangle }{T|g'|^2}=\Gamma_{12}^2
\end{equation*}
analogously
\begin{equation*}
\Gamma_{22}^2=\frac{T\langle g',ig''\rangle-2|g'|^2\langle g,ig'\rangle }{T|g'|^2}=\Gamma_{21}^1
\end{equation*}
On the other hand, $\Gamma_{ii}^j=-\frac{L_{ii,j}}{2L_{jj}}$ for  $i\neq j$, we get
\begin{equation*}
\Gamma_{11}^2=\frac{2|g'|^2\langle g,ig'\rangle -T\langle g',ig''\rangle}{T|g'|^2},\quad\Gamma_{22}^1=\frac{2|g'|^2\langle g,g'\rangle-T\langle g',g''\rangle }{T|g'|^2}.
\end{equation*}
With what has been done so far, we find
\begin{equation*}
 \sum_{k=1}^{2}\Gamma_{11}^k h_{,k} =e^{\langle 1,f\rangle}\langle 1,\xi+f'' \rangle, \quad \sum_{k=1}^{2}\Gamma_{22}^kh_{,k}=-e^{\langle 1,f\rangle}\langle 1,\xi+f'' \rangle, \quad \sum_{k=1}^{2}\Gamma_{12}^kh_{,k}=-e^{\langle 1,f\rangle}\langle i,\xi +f''\rangle,
\end{equation*}
Where $\xi =f' \left(\frac{g''}{g'}-\frac{2}{T}g'\overline{g}\right)-f''$. Using  (\ref{2})  we can write
\begin{equation*}
   V_{11}=\frac{T^2h}{4|g'|^2}[\langle 1,f'\rangle^2-\langle 1,\xi \rangle] + h,\quad V_{22}=\frac{T^2h}{4|g'|^2}[\langle 1,if'\rangle^2+\langle 1,\xi \rangle] + h,\quad V_{12}=\frac{T^2h}{4|g'|^2}[\langle i,\xi-\frac{f'^2}{2} \rangle]
\end{equation*}
Therefore we have (\ref{3}) e (\ref{33}).
\end{Proof}

\begin{teo}
\label{te3}
Let $ \Sigma $ be a Riemann surface and $ X: \Sigma \rightarrow \mathbb{R}^3 $ an immersion such that the Gauss-kronecker curvature is non-zero. \\
Then $ X (\Sigma) $ is a SS-surface if and only if, there are holomorphic functions $ f, g: \mathbb{C} \rightarrow \mathbb{C} $ with $ g' \neq 0 $ , such that $ X (\Sigma) $ is locally parameterized by
\begin{equation*}
\begin{split}
X=\frac{e^{\langle 1,f\rangle}}{2|g'|}(Tg'\bar{f'}-2g\langle g',gf'\rangle , -2\langle g',gf'\rangle) + e^{\langle 1,f\rangle}\frac{(2g,2-T)}{T}
\end{split}
\end{equation*}
With normal application of Gauss $N=\frac{(2g,2-T)}{T}$, $T=1+|g|^2$.
The $ P $ regularity condition is given by
\begin{equation*}
P=detV=\frac{1}{K}\neq 0
\end{equation*}
The coefficients of the first and second fundamental forms of $ X $
they have the following expressions
\begin{equation*}
\begin{split}
E&=\langle X_{,1},X_{,1} \rangle= \frac{T^2h^2}{4|g'|^2}\left((\langle 1,f'\rangle^2-\langle 1,\xi \rangle)^2+\langle i,\xi-\frac{f'^2}{2} \rangle^2\right)+2h^2(\langle 1,f'\rangle^2-\langle 1,\xi \rangle)+\frac{4h^2|g'|^2}{T^2}\\
F&=\langle X_{,1},X_{,2} \rangle=\left(\frac{T^2h^2|f'|^2}{4|g'|^2} +2h^2\right)\langle i,\xi-\frac{f'^2}{2} \rangle\\
   G&=\langle X_{,2},X_{,2} \rangle=\frac{T^2h^2}{4|g'|^2}\left((\langle 1,if'\rangle^2+\langle 1,\xi \rangle)^2+\langle i,\xi-\frac{f'^2}{2} \rangle^2\right)+2h^2(\langle 1,if'\rangle^2+\langle 1,\xi \rangle)+\frac{4h^2|g'|^2}{T^2}\\
e&=\langle X_{,1},N_{,1}\rangle=h(\langle 1,f'\rangle^2-\langle 1,\xi \rangle) + \frac{4h|g'|^2}{T^2}\\
f&=\langle X_{,1},N_{,2} =h\langle i,\xi-\frac{f'^2}{2} \rangle\\
g&=\langle X_{,2},N_{,2} \rangle=h(\langle 1,if'\rangle^2+\langle 1,\xi \rangle) +\frac{4h|g'|^2}{T^2}
\end{split}
\end{equation*}
\end{teo}
\begin{Proof}
Since $ \Sigma $ is a Riemann surface, that is, there is a holomorphic function $ g: \mathbb{C} \rightarrow \mathbb{C}_{\infty} $ and $ g'\neq 0 $ such that $ N:U\subset \mathbb{R}^3 \rightarrow \mathbb{S}^2 $ given by $ N (u) = (Pr ^ {- 1} \circ g) (u) $ e $ Pr^{- 1} $ is the inverse of the stereographic projection.In this way, we have
\begin{equation*}
N(u)=\frac{(2g(u),1-|g(u)|^2)}{1+|g(u)|^2}
\end{equation*}
$X$ is parameterized locally by by equation (\ref{1})
By theorem  \ref{te2} we have  $W=V^{-1}$, let $ \sigma_i $ be the eigenvalues of $ V $ and $ \lambda_i $ the eigenvalues of $ W $ then $ \lambda_i = \frac{1} {\sigma_{i}} $, the same as $ \sigma_ {i} = \frac{1}{\lambda_i} $, like this
\begin{equation}\label{5}
   trV=V_{11}+V_{22}=\frac{-2H}{K}
\end{equation}
On the other hand, of (\ref{2}) we have to
\begin{equation*}
V_{ij}=\frac{1}{L_{jj}}\left(h_{,ij}-\sum_{k=1}^{n}h_{,k}\Gamma_{ij}^k+hL_{ij}\delta_{ij}\right)
\end{equation*}
By the previous lemma we have
\begin{equation}
\begin{split}\label{6}
trV&=V_{11}+V_{22} =\frac{h_{,11}}{L_{11}}-\frac{h_{,1}}{L_{11}}\Gamma_{11}^1-\frac{h_{,2}}{L_{11}}\Gamma_{11}^2+ h+
 \frac{h_{,22}}{L_{22}}-\frac{h_{,1}}{L_{22}}\Gamma_{22}^1-\frac{h_{,2}}{L_{22}}\Gamma_{22}^2+ h\\
 &=\frac{h_{,11}}{L_{11}} + \frac{h_{,22}}{L_{22}} - \frac{h_{,1}}{L_{11}}(\Gamma_{11}^1-\Gamma_{12}^2)-\frac{h_{,2}}{L_{22}}(\Gamma_{22}^2-\Gamma_{21}^1)+ 2h=\frac{h_{,11}}{L_{11}} + \frac{h_{,22}}{L_{22}}+2h\\
 &=\frac{\triangle h}{L_{11}}+2h
\end{split}
\end{equation}

In (\ref{5}) and (\ref{6}) we get
\begin{equation*}
 \left(\frac{h\triangle h-|\nabla h|^2}{L_{11}}\right)K+ 2H\langle X,N\rangle+ (\langle X,X\rangle+\langle X,N\rangle^2)K=0
\end{equation*}
\begin{equation}\label{7}
 \left(\frac{h\triangle h-|\nabla h|^2}{L_{11}}\right)K+ 2H\psi+ (\Lambda+\psi^2)K=0,\quad \forall p\in \Sigma
 \end{equation}
Therefore $ X (\Sigma) $ is a SS-surface if and only if
\begin{equation}\label{lap}
 h \triangle h- | \nabla h |^ 2 = 0
\end{equation}
Let $ h $ be a solution of (\ref{lap}) , without loss of generality we can assume that $ h = e^{\phi} $ where $ \phi: \Sigma \rightarrow \mathbb{R} $ is a differentiable function, in this case,
\begin{equation*}
h_{,1}=e^\phi \phi_{,1}, \quad h_{,2}=e^\phi \phi_{,2},\quad h_{,11}=e^\phi (\phi_{,1}^2+\phi_{,11}),\quad h_{,22}=e^\phi (\phi_{,2}^2+\phi_{,22}).
\end{equation*}
In  (\ref{lap})  we get $\triangle \phi=0$, $\phi$ is a real part of a holomorphic function, such that $\phi=\langle 1,f\rangle$. Therefore $h=e^{\langle 1,f\rangle}$.\\
We consider $L_{ii}=\langle N_{,i},N_{,i}\rangle ,1\leq i\leq 2 $. With $T=1+|g|^2$, we get $T_{,1}=2\langle g,g' \rangle$  e $T_{,2}=2\langle g,ig'\rangle$. So,
\begin{equation*}
   N_{,1}=\frac{2}{T^2}(Tg_{,1}-2g\langle g',g\rangle , -2\langle g',g \rangle),\quad \langle N_{,1},N_{,1} \rangle=\frac{4|g'|^2}{T^2}
\end{equation*}
Where we use the fact that $ g'= g_{,1} $, since $ g $ is holomorphic. Similarly we have
\begin{equation*}
   N_{,2}=\frac{2}{T^2}(Tg_{,2}-2g\langle g,ig'\rangle , -2\langle g,ig' \rangle  ),\quad  \langle N_{,2},N_{,2}\rangle =\frac{4|g'|^2}{T^2}
\end{equation*}
Where we use the equality $ ig'= g_{, 2} $, since $ g $ is holomorphic and using the properties of holomorphic functions. Finally
\begin{equation*}
  \langle N_{,1},N_{,2}\rangle =\frac{4}{T^2} \langle g_{,1},g_{,2} \rangle=0
\end{equation*}
like this
\begin{equation*}
  L_{ij}=\frac{4|g'|^2}{T^2}\delta_{ij}, \hspace{0.5cm} T=1+|g|^2, \hspace{0.5cm} 1\leq i\leq 2
\end{equation*}
Then we can rewrite (\ref{1}) as
\begin{equation*}
   X(q) =\sum_{j=1}^{2} \frac{h_{,j}}{L_{jj}}N_{,j}+hN=\frac{e^{\langle 1,f\rangle}}{2|g'|^2} (Tg'\bar{f'}-2g\langle g',gf'\rangle,-2\langle g',gf'\rangle)+e^{\langle 1,f\rangle}\frac{(2g,2-T)}{T}
\end{equation*}
Using eqaution (\ref{3}) and $ X_{,i}=\sum_{j=1}^{2}V_{ij}N_{,j}$ we obtain
\begin{equation*}
\begin{split}
E&=\langle X_{,1},X_{,1} \rangle=\langle V_{11}N_{,1}+V_{12}N_{,2} ,V_{11}N_{,1}+V_{12}N_{,2} \rangle=V_{11}^2L_{11}+V_{12}^2L_{22}\\
&=\frac{T^2h^2}{4|g'|^2}\left((\langle 1,f'\rangle^2-\langle 1,\xi \rangle)^2+\langle i,\xi-\frac{f'^2}{2} \rangle^2\right)+2h^2(\langle 1,f'\rangle^2-\langle 1,\xi \rangle)+\frac{4h^2|g'|^2}{T^2}
\end{split}
\end{equation*}
Similarly we have
\begin{equation*}
F=\left(\frac{T^2h^2|f'|^2}{4|g'|^2} +2h^2\right)\langle i,\xi-\frac{f'^2}{2} \rangle
\end{equation*}
\begin{equation*}
G=\frac{T^2h^2}{4|g'|^2}\left((\langle 1,if'\rangle^2+\langle 1,\xi \rangle)^2+\langle i,\xi-\frac{f'^2}{2} \rangle^2\right)+2h^2(\langle 1,if'\rangle^2+\langle 1,\xi \rangle)+\frac{4h^2|g'|^2}{T^2}
\end{equation*}
\begin{equation*}
e=h(\langle 1,f'\rangle^2-\langle 1,\xi \rangle) + \frac{4h|g'|^2}{T^2},\quad f=h\langle i,\xi-\frac{f'^2}{2} \rangle,\quad g=h(\langle 1,if'\rangle^2+\langle 1,\xi \rangle) +\frac{4h|g'|^2}{T^2}
\end{equation*}
\end{Proof}
\begin{remark}
Every sphere of radius $r$ at the origin is SS-surface.
\end{remark}

For some holomorphic functions $f$ and $g$ we show some examples of
\begin{equation*}
X=\frac{e^{\langle 1,f\rangle}}{2|g'|}(Tg'\bar{f'}-2g\langle g',gf'\rangle , -2\langle g',gf'\rangle) + e^{\langle 1,f\rangle}\frac{(2g,2-T)}{T}
\end{equation*}

\begin{figure}[htbp]
\begin{center}
\subfigure{\includegraphics[width=32mm]{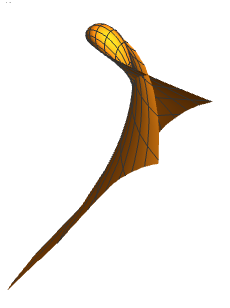}}
\subfigure{\includegraphics[width=32mm]{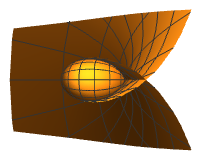}}
\subfigure{\includegraphics[width=32mm]{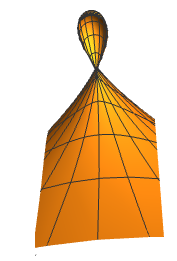}}
\end{center}
\caption{$f(z)=g(z)=z=u_1+iu_2$}
\label{lash3}
\end{figure}

\begin{figure}[htbp]
\begin{center}
\subfigure{\includegraphics[width=45mm]{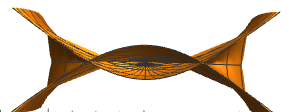}}
\subfigure{\includegraphics[width=29mm]{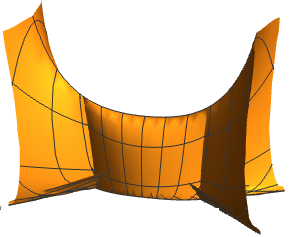}}
\end{center}
\caption{$f(z)=z^2,\hspace{0.2cm} g(z)=z$}
\label{lash3}
\end{figure}

\begin{figure}[htbp]
\begin{center}
\subfigure{\includegraphics[width=32mm]{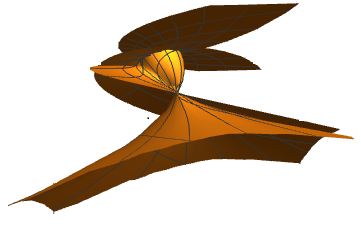}}
\subfigure{\includegraphics[width=28mm]{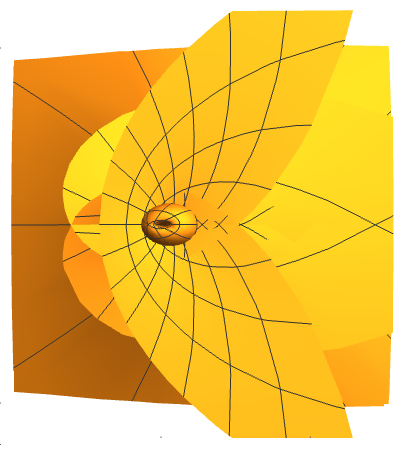}}
\subfigure{\includegraphics[width=32mm]{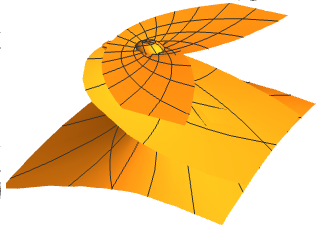}}
\end{center}
\caption{$f(z)=z,\hspace{0.2cm} g(z)=z^3$}
\label{lash3}
\end{figure}

\begin{figure}[htbp]
\begin{center}
\subfigure{\includegraphics[width=32mm]{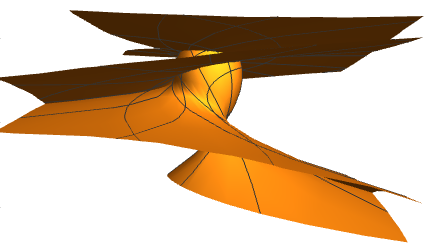}}
\subfigure{\includegraphics[width=28mm]{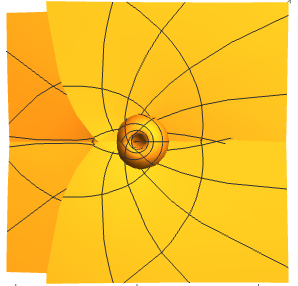}}
\subfigure{\includegraphics[width=32mm]{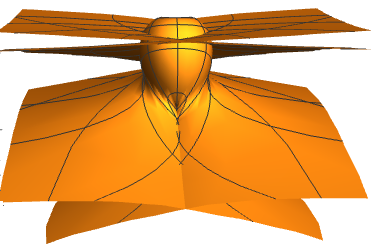}}
\end{center}
\caption{$f(z)=z,\hspace{0.2cm} g(z)=z^4$}
\label{lash3}
\end{figure}
\newpage
The following theorem characterizes the rotating SS-surfaces.
\begin{teo}
Let $ \Sigma $ be a connected SS-surface. since $ \Sigma $ is of rotation if, and only if, there are constants $ a, b \in \mathbb{R} $, such that $ \Sigma $ can be locally parameterized by
\begin{equation}\label{3.5}
X_{a,b}(u_1,u_2)=\left(M(u_1)\cos(u_2),M(u_1)\sin(u_2),N(u_1)\right)
\end{equation}
Where
\begin{equation}\label{3.51}
M(u_1)=e^{au_1+b} \left[\frac{ a(1-e^{2u_1})}{2}+\frac{2e^{u_1}}{1+e^{2u_1}}\right],\quad N(u_1)=e^{au_1+b} \left[\frac{1-e^{2u_1}}{1+e^{2u_1}}-ae^{u_1}\right]
\end{equation}
\end{teo}

\begin{Proof}
By theorem \ref{te3} we have $\Sigma$ is locally parameterized by
\begin{equation*}
X=\frac{e^{\langle 1,f\rangle}}{2|g'|}(Tg'\bar{f'}-2g\langle g',gf'\rangle , -2\langle g',gf'\rangle) + e^{\langle 1,f\rangle}\frac{(2g,2-T)}{T}
\end{equation*}
With normal application of Gauss $ N = \frac{(2g, 2-T)}{T} $, $ T = 1 + | g |^ 2 $ and where $ f, g $ are holomorphic functions. that $ \Sigma $ is of rotation if, and only if, $ g (w) = w $ and $ h (w) = J (| w |^ 2) $, $ w \in \mathbb{C} $, for some $ J $ differentiable function.\\
Changing parameters
\begin{equation*}
w=e^z, \hspace{2cm} z=u_1+iu_2\in \mathbb{C}
\end{equation*}
We will have $ g (z) = z $ and $ h (z) = J(e^{2u_1}) $, consequently, $ h,_2 = 0 $, remembering that $ h = e^{\langle 1, f \rangle} $ e $ h,_2 = e^{\langle 1, f \rangle} \langle 1, if'\rangle = 0 $, then $ \langle 1, if' \rangle = 0 $, but $ \langle 1, if'\rangle = \langle 1, f \rangle,_2 = 0 $ and so on by the Cauchy-Riemann equations, $ \langle i, f \rangle, _1 = 0 $. \\
On the other hand, the real part of a holomorphic function is harmonic
\begin{equation*}
  \langle 1, f \rangle,_{11}+\langle 1, f \rangle,_{22}=0
\end{equation*}
As $\langle 1,f\rangle,_2=0$ then $\langle 1,f\rangle,_{22}=0$ we have $\langle 1,f\rangle,_{11}=0$, then,  $\langle 1,f\rangle,_{1}=a$ with $a\in \mathbb{R}$, Which shows that $\langle 1, f \rangle=au_1+b$, $a,b\in \mathbb{R}$.\\
That way, we have to
\begin{equation*}
f(z)=az+z_0,\hspace{.2cm}g(z)=e^z, \hspace{.2cm} h(z)=e^{au_1+b},\hspace{.2cm}z=u_1+iu_2,z_0=b+ic\in \mathbb{C}
\end{equation*}
So that $f'(z)=a$, $g'(z)=g(z)$ e $T=1+|g|^2=1+|e^{2u_1}(\cos(u _2)+i \sin(u_2))|^2=1+e^{2u_1}$.In these conditions we have to
\begin{equation*}
\begin{split}
   X &=\frac{e^{\langle 1,f\rangle}}{2|g'|}(Tg'\bar{f'}-2g\langle g',gf'\rangle , -2\langle g',gf'\rangle) + e^{\langle 1,f\rangle}\frac{(2g,2-T)}{T}
   \\   &=e^{au_1+b}\left(\left(\frac{a(1-e^{2u_1})}{2}+\frac{2e^{u_1}}{1+e^{2u_1}}\right)(\cos(u_2)+i\sin(u_2)),\frac{1-e^{2u_1}}{1+e^{2u_1}}-ae^{u_1}\right)
\end{split}
\end{equation*}
Therefore  $\Sigma$ can be parameterized locally by (\ref{3.5}) where (\ref{3.51})
\end{Proof}
The following are examples of the previous theorem for different types of parameters $a$ e $b$
\begin{figure}[htbp]
\begin{center}
\subfigure{\includegraphics[width=32mm]{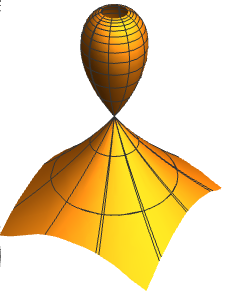}}
\subfigure{\includegraphics[width=32mm]{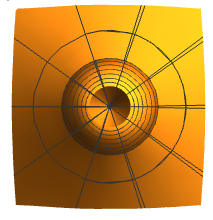}}
\end{center}
\caption{$a=1,\hspace{0.2cm} b=0$}
\label{lash3}
\end{figure}

\begin{figure}[htbp]
\begin{center}
\subfigure{\includegraphics[width=32mm]{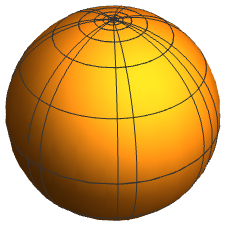}}
\end{center}
\caption{$a=0,\hspace{0.2cm} b\in \mathbb{R}$}
\label{lash3}
\end{figure}

\begin{figure}[htbp]
\begin{center}
\subfigure{\includegraphics[width=32mm]{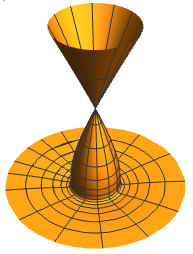}}
\subfigure{\includegraphics[width=32mm]{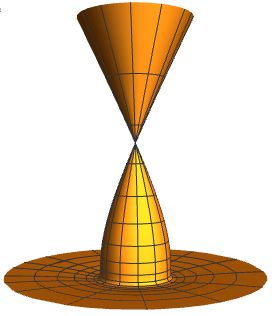}}
\end{center}
\caption{$a=-1,\hspace{0.2cm} b=1$}
\label{lash3}
\end{figure}
\newpage
\begin{teo}
Let $ \Sigma $ be a compact, connected SS-surface then $ \Sigma $ is a sphere.
\end{teo}
\begin{Proof}
There is a $ E $ sphere, such that $ \Sigma $ is contained within the sphere, and a point $ p \in E \cap \Sigma $, such that $ T_pE = T_p\Sigma $.
Let $ h_1, h_2: U \subset \mathbb{R}^2 \rightarrow \mathbb{R} $, where $ h_1 $ the support function of $ E $ in this case $ h_1 $ be constant and $ h_2 $ function support of $ \Sigma $. In addition $ h_1 (p) = h_2 (p) $. We know that $ h_2 (q) \leq |q| $, for every $ q \in \Sigma $. So
\begin{equation*}
h_2(q)\leq |q| \leq |p|=h_1(p) \quad \text{for all} \quad q\in \Sigma
\end{equation*}
\begin{equation*}
e^{\mu_2}\leq e^{\mu_1} \Rightarrow \mu_2\leq \mu_1
\end{equation*}
Where $ \mu_1 $ and $ \mu_2 $ are harmonics and $ \mu_2(p) = \mu_1(p) $, just by the principle of the maximum, $ \mu_2 = \mu_1 $ in $ U $. Hence $ \mu_2 $ is constant, so $ h_2 $ is constant and therefore $ \Sigma $ is a sphere.
\end{Proof}

\newpage

\renewcommand{\refname}{Bibliografía}

\end{document}